\documentclass{article}

\usepackage{amsmath}
\usepackage{amssymb}
\usepackage{amsopn}
\usepackage{epsfig}
\usepackage{amsfonts}
\usepackage{latexsym}
\usepackage{graphicx}

\usepackage{color}
\definecolor{darkgreen}{rgb}{0,.4,0.2}
\definecolor{darkagenta}{rgb}{.5,0,.5}
\definecolor{darkred}{rgb}{0.85,0,0}
\definecolor{darkblue}{rgb}{0,0,.6}
\definecolor{lightgray}{gray}{.95}

\newcommand{\red}[1]{{\color{red}#1}}
\newcommand{\blue}[1]{{\color{blue}#1}}
\newcommand{\darkred}[1]{{\color{darkred}#1}}

\newcommand{\darkgreen }[1]{{\color{darkgreen}#1}}
\newcommand{\darkagenta}[1]{{\color{darkagenta}#1}}

\newcommand{\Sh}{{\rm S}}
\newcommand{\Lan}{{\mathcal{L}}}
\newcommand{\Id}{ {\rm Id}}
\newcommand{\N}{{\mathbb{N}}}
\newcommand{\Z}{{\mathbb{Z}}}

\newtheorem{theorem}{Theorem}[section]

\newtheorem{lemma}{Lemma}[section]
\newtheorem{proposition}{Proposition}[section]

\date{\today}

\begin{document}

\begin{center}
{\bf \Large The Frobenius problem for  homomorphic\\[.3cm] embeddings of languages into the integers}

\bigskip

{\bf \large Michel Dekking}

\bigskip

{\footnotesize DIAM, Delft University of Technology, December 12, 2017 }

\end{center}

\section{Introduction}

\medskip

\qquad The Frobenius problem is also known as the `coin problem'. Since the value of a coin can only be positive, we will consider exclusively embeddings into the natural numbers $\mathbb{N}=\{1,2,3,\dots\}$. Let $\Lan$ be a language, i.e., a sub-semigroup of the free semigroup generated by a finite alphabet under the concatenation operation.

\qquad A homomorphism
of $\Lan$ into the natural numbers is a map $\Sh:\Lan\rightarrow \N$  satisfying $$\Sh(vw)=\Sh(v)+\Sh(w), \quad {\rm for\: all }\: v,w \in \Lan.$$

The two main questions to be asked about the image set $\Sh(\Lan)$ are

\smallskip

(Q1) Is the complement $\N\setminus \Sh(\Lan)$  finite or infinite?

\smallskip

(Q2) If the complement of $\Sh(\Lan)$ is finite, then what is the largest element in this set?

\smallskip

These two questions are known as the Frobenius problem in the special case that $\Lan$ is the full language consisting of \emph{all} words over a finite alphabet.
In this case they have been posed as a problem (with solution) for an alphabet $\{a,b\}$ of cardinality 2 by James Joseph Sylvester in 1884 \cite{Sylvester}: $\N\setminus\Sh(\Lan)$ is finite, and its largest element is
$$\Sh(a)\Sh(b)-\Sh(a)-\Sh(b). $$
In this paper we will also restrict ourselves to the two symbol case: alphabet $\{a,b\}$.

\qquad In Section \ref{sec:gold} we prove that for the golden mean language (``no $bb$") the set
$\N\setminus\Sh(\Lan)$ is finite, with largest element
$$\Sh(a)^2+\Sh(a)\Sh(b)-3\Sh(a)-\Sh(b). $$
Our main interest is however not in sofic languages\footnote{Languages defined by the labelling of infinite paths of an automaton.}, but in languages with low complexity, where the complement of $\Sh(\Lan)$ can be infinite.

\qquad In Section \ref{sec:Sturm} we analyse the case of Sturmian languages, and show that for the Fibonacci language a $0\!-\!\infty$ law holds: either the complement is empty or it has infinite cardinality.

\qquad  In Section \ref{sec:TM} we show that for any homomorphism $\Sh$ the image of the Thue-Morse language will consist of a union of 5 arithmetic sequences.

\qquad  In Section \ref{sec:twodim} we consider two-dimensional embeddings, which behave quite differently.

\medskip

\qquad We usually suppose that $\gcd(\Sh(a),\Sh(b))=1$. First of all this is not a big loss since automatically the complement will have infinite cardinality in this case. Secondly, if $r$ divides both $\Sh_1(a)$ and $\Sh_1(b)$ for some homomorphism $\Sh_1$, then
$$\Sh_1(\Lan^n)=r^n\Sh_2(\Lan^n), \quad {\rm for}\: n=1,2,\dots,\: {\rm where}\; \Sh_2(a)=\frac{\Sh_1(a)}{r},\, \Sh_2(b)=\frac{\Sh_1(b)}{r}. $$

\medskip

\qquad  Our work is related to the work on  \emph{abelian complexity}, see, e.g., \cite{abeliancomp}, \cite{RSZ}, \cite{Karhumaki-S-K}.
See Lemma \ref{lem:compl} for such a connection.

\medskip

\qquad  Our work is also related to the notion of \emph{additive complexity}, see \cite{Sahas} and \cite{Ardal}. The \emph{additive complexity} of an infinite word $w$ over a finite set of
integers (see \cite{Ardal}) is the function $n\rightarrow\phi^+(w, n)$ that counts the number of distinct sums obtained
by summing $n$ consecutive symbols of $w$.
In general we write $\Lan^n$ for the set of words of length $n$ in a language $\Lan$. Let $\Lan_w$ be the language of all words occurring in the infinite word $w$. Then the additive complexity is  $\phi^+(w, n)={\rm Card} \{\Sh(u): u\in \Lan_w^n\}$, where $\Sh$ is the identity map on the alphabet of $w$.

\medskip

\qquad  We finally mention that homomorphisms $\Sh$ from a language to the natural numbers already occur in the 1972 paper \cite[Section 6]{Car-Sco-Hog} in the context of the Fibonacci language, where they are called \emph{weights}.

\section{Homomorphic images of the golden-mean language}\label{sec:gold}

\qquad The golden mean language is the language $\Lan_{ \rm GM}$ consisting of all words over $\{a,b\}$ in which $bb$ does not occur as a subword.
Now if $\Sh$ satisfies $\Sh(a)=1$ or  $\Sh(b)=1$, then it is easily seen that $\Sh(\Lan_{ \rm GM})=\N$,  so for these homomorphisms the golden mean and the full language both map to $\mathbb{N}$. One could say they both have Frobenius number $0$. In general however, the Frobenius number will increase substantially. If we take $\Sh$ defined by
$$\Sh(a)=100,\; \Sh(b)=3,$$
then the Frobenius number of the full language under $\Sh$ is $300-100-3=197$, and the Frobenius number of $\Sh(\Lan_{ \rm GM})$ is equal to 9997.  For arbitrary homomorphisms the solution of the Frobenius problem for the golden mean language is given by  the following, where we write $\Sh_a:=\Sh(a),\, \Sh_b:=\Sh(b)$.

 \begin{theorem}\label{th:Goldmean} Let  $\Sh: \Lan_{ \rm GM}\rightarrow \N$ be a homomorphism.  Suppose \,$\gcd(\Sh_a,\Sh_b)=1$, and both $\Sh_a>1$ and  $\Sh_b>1$. Then the Frobenius number of $\Sh(\Lan_{ \rm GM})$ is equal to
 $$\max\: \N\setminus \Sh(\Lan_{ \rm GM})=  \Sh_a(\Sh_a-3) +\Sh_b(\Sh_a-1).$$
\end{theorem}

\emph{Proof:} Let an $\Sh_a$-$point$ be defined as a multiple $n\Sh_a$, $n=0,1,\dots$, and an $\Sh_a$-$interval$ as the set of numbers between two consecutive $\Sh_a$-points. We also consider $\Sh_b$-$chains$, defined for  $n\ge0$ by
$$C(n)=\{n\Sh_a+\Sh_b,\, n\Sh_a+2\Sh_b,\dots, n\Sh_a+(n+1)\Sh_b\}.$$
Note that the union of the  $\Sh_a$-points and the $\Sh_b$-chains will give $\Lan_{ \rm GM}$.

The key observation is that the $\Sh_b$-chain $C(\Sh_a-2)$ has
$\Sh_a-1$ elements, which are all different modulo $\Sh_a$. This is a consequence of $\gcd(\Sh_a,\Sh_b)=1$. It follows that the $\Sh_b$-chains fill in more and more points of the $\Sh_a$-intervals. The last point to be filled in is modulo $\Sh_a$ equal to $\Sh_a-\Sh_b$, produced by the last element of the chain $C(\Sh_a-2)$.

This is the number
$$P:=(\Sh_a-2)\Sh_a+(\Sh_a-1)\Sh_b.$$
But then the largest number in the complement of $\Lan_{ \rm GM}$ is $P-\Sh_a$, which is the number as claimed in the theorem. In this argument we used that if  a point in an $\Sh_a$-interval is filled in, then the corresponding points modulo $\Sh_a$ in all later intervals will also be filled in, simply because the later chains will be extensions of the earlier ones. \hfill$\Box$

\begin{figure}[!h]
\begin{center}
\includegraphics[height=3.5cm]{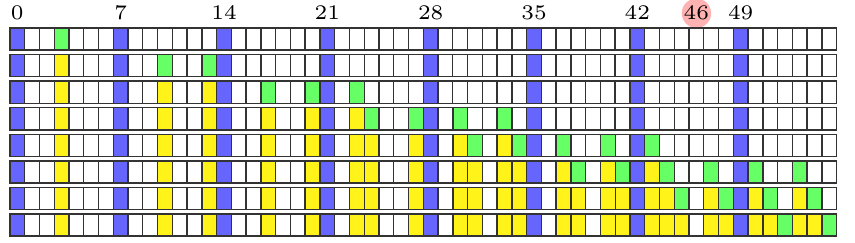}
\vspace*{-.5cm}
\end{center}
\caption{\small Example with $\Sh(a) = 7, \Sh(b) = 3$: row $n$ shows the $\Sh_a$-points in blue, and the $\Sh_b$-chain $C(n-1)$ in yellow and green, for $n = 1,\dots,8$ (truncated at 56).}\label{gold}
\end{figure}

 \section{Sturmian languages}\label{sec:Sturm}

 \qquad Sturmian words are infinite words over a two letter alphabet that haven exactly $n+1$ subwords for each $n=1,2,\dots$. We call the collection of these subwords a Sturmian language. There is a surprising characterization of Sturmian words: $s$ is Sturmian if and only if $s$ is irrational \emph{mechanical}, which means that there exists an irrational number $\alpha\in (0,1)$ and a number $\rho$ such that $s=s_{\alpha,\rho}$, or $s=s'_{\alpha,\rho}$, where
 $$s_{\alpha,\rho}=\big([ (n+1)\,\alpha +\rho]-[ n\,\alpha +\rho]\big)_{n\ge 0},\;
    s'_{\alpha,\rho}=\big(\lceil (n+1)\,\alpha +\rho\rceil-\lceil n\,\alpha +\rho\rceil\big)_{n\ge 0}. $$
 See, e.g., \cite[Prop.~2.1.13]{Lothaire}.
 Because of this representation, we will use the alphabet $\{0,1\}$ instead of $\{a,b\}$ in this section.

 Of special interest are the Sturmian words $s_\alpha:=s_{\alpha,0}$  and $s'_\alpha:=s'_{\alpha,0}$ of intercept 0. These have the property that they only differ in the first element:
 $$s_\alpha=0\,c_\alpha, \qquad s'_\alpha=1\,c_\alpha.$$
  Here $c_\alpha:=s_{\alpha,\alpha}$ is called the \emph{characteristic word} of $\alpha$.
  For $n\ge 0$ we have
 $$c_\alpha(n)=s_{\alpha,\alpha}(n)=[ (n+1)\,\alpha +\alpha]-[ n\,\alpha +\alpha]
 =[ (n+2)\,\alpha ]-[ (n+1)\,\alpha ].$$

The words $s_\alpha, s'_\alpha$ and $c_\alpha$ generate the same language (\cite[Prop.2.1.18]{Lothaire} ), which we denote $\Lan_\alpha$. Recall that $\Lan^n_\alpha$ is the set of words of length $n$ in $\Lan_\alpha$.

 \begin{lemma}\label{lem:compl} Let $\Lan_\alpha$ be a Sturmian language, and  let $\Sh$ be a homomorphism with $\Sh(0)\ne\Sh(1)$. Then ${\rm Card} \,\Sh(\Lan^n_\alpha)=2$ for all $n\ge 1$.  \end{lemma}

 \emph{Proof:} This follows directly from the fact (\cite[Th.2.1.5]{Lothaire}) that Sturmian words are \emph{balanced}, i.e., any two words of the same length can at most differ 1 in their number of ones.\hfill$\Box$

 \medskip

 A sequence $([n\alpha])$, where $[.]$ denotes integer part, is called a \emph{Beatty sequence} if $\alpha>1$, and a \emph{slow Beatty sequence} if  $0<\alpha<1$
 (terminology from \cite{KimStol}).

  \begin{theorem}\label{th:Sturm} Let  $\alpha$ be an irrational number from $(0,1)$. Let $\Lan_\alpha$ be the Sturmian language generated by $\alpha$, and let  $(q_n)_{n\ge 0}$ be the slow Beatty sequence defined by $$q_n=[(n+1)\alpha].$$ Let $\Sh: \Lan_\alpha\rightarrow \N$ be a homomorphism. Define $\Sh_0=\Sh(0), \Sh_1=\Sh(1)$. Then
 $$ \Sh(\Lan_\alpha) = \{(\Sh_1-\Sh_0)q_n+n\Sh_0+\Sh_0: \; n=0,\dots\} \cup \{(\Sh_1-\Sh_0)q_n+n\Sh_0+\Sh_1: \; n=0,\dots\}.$$
\end{theorem}

\emph{Proof:}  If $\Sh_0=\Sh_1$ then this is certainly true, so suppose $\Sh_0\ne\Sh_1$ in the sequel. We denote $c_\alpha[i,j]:=c_\alpha(i)\dots c_\alpha(j)$ for integers $0\le i<j$. Let $N_{\ell}(w)$  denote the number of occurrences of the letter $\ell$ in a word $w$ for $\ell=0,1$.
Then
$$N_1(c_\alpha[0, n-1])=\sum_{k=0}^{n-1} c_\alpha(k)=[(n+1)\alpha]-[\alpha]=q_n, \quad N_0(c_\alpha[0,n-1])=n-q_n.$$
Of course all words $c_\alpha[0,n-1]$ are in the Sturmian language $\Lan_\alpha$,  but  $\Lan_\alpha$ also contains the words $0c_\alpha[0,n-1]$ and $1c_\alpha[0,n-1]$. It thus follows from Lemma \ref{lem:compl} that $\Sh(\Lan_\alpha)$ is given by the union of all images  $\Sh(0c_\alpha[0,n-1])$ and $\Sh(1c_\alpha[0,n-1])$. Since $$\Sh(0c_\alpha[0,n-1])=\Sh_0+(n-q_n)\Sh_0+q_n\Sh_1=(\Sh_1-\Sh_0)q_n+n\Sh_0+\Sh_0,$$  the result follows.\hfill$\Box$

\medskip

\subsection{The Fibonacci language }

\qquad Let $\Phi=(\sqrt{5}+1)/2=1.61803\dots$ be the  golden mean, and let $\alpha:=2-\Phi$. We have
$$c_\alpha=([(n+1)\alpha]-[n\alpha])_{n\ge 1}= 0, 1, 0, 0, 1, 0, 1, 0, 0, 1, 0, 0, 1, 0, 1, 0, 0,\dots,$$
the infinite Fibonacci word. We write $\Lan_{\rm F}:=\Lan_\alpha$.

 \begin{theorem}\label{th:Fib} Let  $\Sh: \Lan_{\rm F}\rightarrow \N$ be a homomorphism.  Then
 $$ \Sh(\Lan_{\rm F})=\big((\Sh_0-\Sh_1)[n\Phi]+(2\Sh_1-\Sh_0)n+\Sh_0-\Sh_1\big)_{n\ge 1} \,\cup\, \big((\Sh_0-\Sh_1)[n\Phi]+(2\Sh_1-\Sh_0)n\big)_{n\ge 1}.$$
\end{theorem}


 \emph{ Proof:} This is a corollary to Theorem \ref{th:Sturm}, using $[-x]=-[x]-1$  for non-integer $x$:
 \begin{eqnarray*}
 (\Sh_1-\Sh_0)q_{n-1}+n\Sh_0 &=& (\Sh_1-\Sh_0)[n\alpha]+n\Sh_0= (\Sh_1-\Sh_0)[n(2-\Phi)]+n\Sh_0\\
  &=& 2(\Sh_1-\Sh_0)n+(\Sh_1-\Sh_0)[-n\Phi]+n\Sh_0\\
  &=& (2\Sh_1-\Sh_0)n+(\Sh_1-\Sh_0)(-[n\Phi]-1)\\
  &=& (\Sh_0-\Sh_1)[n\Phi]+(2\Sh_1-\Sh_0)n+\Sh_0-\Sh_1.\hspace*{4cm}\Box
\end{eqnarray*}

\medskip

\begin{lemma}\label{lem:empty}  For $\Sh(0)=1,\, \Sh(1)\le 3$ or $\Sh(0)=2,\, \Sh(1)=1$ one has $\Sh(\Lan_{\rm F})=\N$.
\end{lemma}

\emph{Proof:}
Take  $(\Sh_0,\Sh_1)=$ (1,1). Then obviously $\Sh(\Lan_{\rm F})=\N$.

Take  $(\Sh_0,\Sh_1)=$ (2,1). Then $\Sh(\Lan_{\rm F})=\N$, since by Theorem~\ref{th:Fib} $\Sh(\Lan_{\rm F})$ is the union of $([n\Phi])$ and $([n\Phi]+1)$, where the difference of two consecutive  terms in $([n\Phi])$ is never more than 2.

Take  $(\Sh_0,\Sh_1)=$ (1,2). Then $\Sh(\Lan_{\rm F})=\N$, since $\Sh(\Lan_{\rm F})$ is the union of $([n(3-\Phi)])$ and $([n(3-\Phi)])+1)$, where the
difference of two consecutive  terms in $([n(3-\Phi)])$ is never more than 2.

Take  $(\Sh_0,\Sh_1)=$ (1,3). This case is more complicated. Let $u:=(-2[n\Phi]+5n-2)_{n\ge 1}$, and $v:=u+2$. Then according to Theorem \ref{th:Fib}, the union of the sets determined by $u$ and $v$ is $\Sh(\Lan_{\rm F})$. Let $\Delta u$ be the difference sequence defined by $\Delta u_n= u_{n+1}-u_n$ for $n \ge 0$. It is easy to see that the difference sequences
$\Delta v$ and $\Delta u$ are both equal to the Fibonacci sequence $1,3,1,1,3,1,\dots$ on the alphabet $\{1,3\}$ (cf.~\cite{AllDekk}). We claim that if two consecutive numbers $m , m+1$ are missing in $u$, then these two do appear in $v$, implying that $\Sh(\Lan_{\rm F})=\N$. Indeed the two missing numbers are characterized by $u_{n+1}-u_n=3$ for some $n$, and the missing numbers are $m=u_n+1$ and $u_n+2$. The second number appears in $v$, simply because $v=u+2$.
The first number appears because $u_{n+1}-u_n=3$ implies $u_n-u_{n-1}=1$ (no 33 in the 1-3-Fibonacci sequence), and so $v_{n-1}=v_n-1=u_n+1$.
\hfill$\Box$

\medskip

We define  $\mathcal{E}:=\{(1,1), (1,2), (1,3), (2,1)\}$.

\medskip

 \begin{theorem}\label{th:card} Let  $\Sh: \Lan_{\rm F}\rightarrow \N$ be a homomorphism. Then
 $\N\setminus \Sh(\Lan_{\rm F})$ has infinite cardinality, \emph{unless} $(\Sh(0),\Sh(1))\in \mathcal{E}$, in which case the complement is empty.
\end{theorem}

\emph{Proof:}  According to Lemma \ref{lem:empty} the complement of $\Sh(\Lan_{\rm F})$ is empty for $(\Sh_0,\Sh_1)\in \mathcal{E}$.

The density of the  set $\Sh(\Lan_{\rm F})$ in the natural numbers exists, and equals
$$\delta:=\frac2{(\Sh_0-\Sh_1)\Phi+\,2\Sh_1-\Sh_0}.$$
The theorem will be proved if we show that $\delta<1$ for $(\Sh_0,\Sh_1)$ not in $\mathcal{E}$. First we note that the denominator of $\delta$ is positive:
$$ (\Sh_0-\Sh_1)(\Phi-1)+\Sh_1>-\Sh_1(\Phi-1)+\Sh_1=\Sh_1(2-\Phi)>0,$$
where we used that $1<\Phi<2$. We now have
$$\delta<1\;\Leftrightarrow\;  (\Sh_0-\Sh_1)\Phi+\,2\Sh_1-\Sh_0 >2   \;\Leftrightarrow\; (\Sh_0-\Sh_1)\Phi > \Sh_0-\Sh_1 +2-\Sh_1.$$
If $\Sh_0>\Sh_1$, this is satisfied, since under this condition $(2-\Sh_1)/(\Sh_0-\Sh_1)\le 0$, unless $(\Sh_0,\Sh_1)=(2,1)\in \mathcal{E}$.
If $\Sh_0<\Sh_1$, we have to see that $\Phi < 1 +(2-\Sh_1)/(\Sh_0-\Sh_1)$. This holds for $\Sh_0\ge2$, since then $(2-\Sh_1)/(\Sh_0-\Sh_1)\ge 1$.
If $\Sh_0=1$, then this does not hold for $\Sh_1=1,2,3$, i.e., for pairs from $\mathcal{E}$,  but it will hold for all $\Sh_1\ge 4$.\hfill$\Box$

 \medskip

 \qquad For particular values of $\Sh(0)$ and $\Sh(1)$ the complement  of the embedding of the language has a nice structure, as it can be expressed in
 the classical Beatty sequences $A(n)=[ n\Phi ]$ for $n\ge 1$, and  $ B(n)=[ n\Phi^2]$ for $n\ge 1$.  The sequences $A$ and $B$ are called the \emph{lower Wythoff sequence} and \emph{upper Wythoff sequence}; they are extremely well-studied.

\medskip

 {\bf Example 1}. Let $\Sh$ be given by  $\Sh(0)=3$ and $\Sh(1)=2$. In the following we use the notation $pX+qY+r= (pX(n)+qY(n)+r)_{n\ge 1}$ for real numbers $p,q,r$ and functions $X,Y:\N\rightarrow\N$. Then

 \qquad $\Sh(\Lan_{\rm F})= B(\N) \:\cup\: B\!+\!1\,(\N), \quad \N\!\setminus\Sh(\Lan_{\rm F})= \{1,4, 9, 12, \dots\}= 2A+\Id+1\,(\N \cup \{0\}).$

 The first statement follows directly from Theorem \ref{th:Fib}. The second statement follows in a number of steps from the fact that $A$ and $B$ form a Beatty pair: $A(\N)\cap B(\N)=\emptyset$, and $A(\N)\cup B(\N)=\N$. This implies that $A(A(\N))\cup A(B(\N))\cup B(\N)=\N$, where the three sets are disjoint. But
 $AA=B-1$ (see, e.g., Formula (3.2) in \cite{Car-Sco-Hog}). Adding 1 to all three sequences it follows that

 \qquad $B(\N) \:\cup\: B+1\,(\N) \:\cup\:  AB+1\,(\N) = \N\setminus \{1\}.$

 Moreover, according to \cite[Formula (3.5)]{Car-Sco-Hog} one has  $ AB=A+B=2A+\Id$.

 But then the three sequences
  $([ n\Phi ]+n)_{n\ge 1},\, ([ n\Phi ]+n+1)_{n\ge 1},\,  (2[ n\Phi ]+n+1)_{n\ge 0}$, form a  complementary triple, i.e., as sets they are disjoint, and their union is $\N$.

  A similar result holds for\footnote{In these two cases  $\N\setminus\Sh(\Lan_{\rm F})$ is given by sequences  A276885, respectively  A276886 in OEIS~(\cite{OEIS-Fib}). It is easily seen that the definitions of these sequences in OEIS are equivalent to the way in which we obtain them.}
   $\Sh(0)=4,\, \Sh(1)=3$.

  \medskip

  {\bf Example 2}. Let $\Sh$ be given by   $\Sh(0)=3$ and $\Sh(1)=1$, then by Theorem \ref{th:Fib}

 \qquad $\Sh(\Lan_{\rm F})= 2A-\Id\,(\N) \:\cup\: 2A-\Id+2\,(\N)$.

 It  is proved in \cite{AllDekk} that

 \qquad$\N\setminus\Sh(\Lan_{\rm F})= \{2, 9, 20, 27, 38, 49, \dots\}= 4A+3\Id+2\,(\mathbb{N} \cup \{0\}),$

 and that the three sequences
  $(2[ n\Phi ]-n)_{n\ge 1},\, (2[ n\Phi ]-n+2)_{n\ge 1},\,  (4[ n\Phi ]+3n+2)_{n\ge 0}$, form a  complementary triple.

\section{The Thue-Morse language}\label{sec:TM}

\qquad Let $\theta$ given by $\theta(a)=ab,\, \theta(b)=ba$ be the Thue-Morse morphism. Let $\Lan_{\rm TM}$ be the language generated by this morphism.

Let $R_{r,s}=\{s,r+s,2r+s,\dots\}$ be the set determined by the arithmetic sequence with terms $rn+s$ for $n=0,1\dots$.

\smallskip

 \begin{theorem}\label{th:TM} Let  $\Sh: \Lan_{\rm TM}\rightarrow \N$ be a homomorphism. Define $p=\Sh(0), q=\Sh(1)$. Then
 $$ \Sh(\Lan_{\rm TM})=R_{p+q,0}\cup R_{p+q,p}\cup R_{p+q,q}\cup R_{p+q,2p}\cup R_{p+q,2q}.$$
\end{theorem}

\emph{Proof:} Let $\Lan_{\rm TM}^n$ be the set of words of length $n$ in the Thue-Morse language. Put $r=\Sh(ab)=p+q.$
It is clear (and for $p=0, q=1$ observed also in \cite{RSZ}) that since the Thue-Morse word is a non-periodic concatenation of $ab$ and $ba$ that for $n=1,2,\dots$
$$\Sh(\Lan_{\rm TM}^{2n})=\{rn, rn+q-p, rn+p-q\},\quad \Sh(\Lan_{\rm TM}^{2n-1})=\{rn+p, rn+q\}.$$
This implies the statement of the theorem.\hfill$\Box$

\bigskip

 \begin{theorem}\label{th:cardTM} Let  $\Sh: \Lan_{\rm TM}\rightarrow \N$ be a homomorphism. Then
 $\N\setminus \Sh(\Lan_{\rm F})$ has infinite cardinality if and only if\, $\Sh(a)+\Sh(b)\ge 6$. For $\Sh(a)+\Sh(b) < 6$,  the complement is either empty or a singleton.
\end{theorem}

\emph{Proof:} This follows directly form Theorem \ref{th:TM}. If $\Sh(a)+\Sh(b)\ge 6$, then the density of $\N\setminus\Sh(\Lan_{\rm TM})$, is at least 1/6, so the set has infinite cardinality. The results for $\Sh(a)+\Sh(b) < 6$ follow also directly from the previous theorem. \hfill$\Box$

\medskip

{\bf Remark} Let $\sigma$ given by $\sigma(a)=ab, \sigma(b)=aa$ be the period-doubling or Toeplitz morphism.
The difficulty---see \cite[Lemma 6]{Karhumaki-S-K}---of determining the abelian complexity of the period-doubling morphism already indicates that solving the Frobenius problem for the period-doubling language will be much more involved than for the Thue-Morse language.

\section{Two dimensional embeddings}\label{sec:twodim}

\qquad Here we consider homomorphisms $\Sh: \Lan \rightarrow \N\times \N$ and  $\Sh: \Lan \rightarrow \Z\times\Z$. The situation changes drastically for this `double-coin' problem.

 \begin{proposition}\label{prop:2D} Let $\Lan$ be a language on the alphabet $\{a,b\}$, and let $\Sh: \Lan\rightarrow \N\times \N$ be a homomorphism. Then  $\N\times\N\,\setminus\, \Sh(\Lan)$ has infinite cardinality for all pairs $\{\Sh(a), \Sh(b)\}$ which are not equal to the pair $\{{ \rm (0,1),(1,0)}\}$.
\end{proposition}

\emph{Proof:} It suffices to prove this for the full language $\Lan_{\rm full}$.
The image under $\Sh$ is an integer lattice, with a complement of infinite cardinality,
unless $\Sh(a)$ and $\Sh(b)$ are the unit vectors. \hfill$\Box$

\medskip

\qquad We learn from this that the alphabet is `too small', and that we should rather consider embeddings in $\Z\times\Z$ instead of $\N\times \N$.
We focus again on low complexity languages, in particular on those generated by a primitive morphism $\varphi$ on an alphabet $A$. Such a morphism
has a language $\Lan_\varphi$ associated to it, where each word $w\in \Lan_\phi$ has a measure $\mu_\varphi(w)$.
For a given homomorphism $\Sh:\Lan_\varphi\rightarrow \Z\times\Z$ we call the average
$$\Delta_\varphi(\Sh):=\sum_{a\in A} \mu_\varphi(a)\Sh(a)\vspace*{-.5cm}$$
the \emph{drift} of $\Sh$.

 \begin{proposition}\label{prop:drift} Let $\Lan_\varphi$ be a language generated by primitive morphism on an alphabet $A$, and let $\Sh: \Lan_\varphi\rightarrow \Z\times \Z$ be a homomorphism. Then  $\Z\times\Z\,\setminus\, \Sh(\Lan)$ has infinite cardinality if $\Delta_\varphi(\Sh)\ne (0,0)$.
\end{proposition}

\emph{Proof:} It is well-known that the measure $\mu_\varphi$ is strictly ergodic. Because of  this, we have for words  $w$ from $\Lan_\varphi$, where $|w|$ denotes the length of $w$,
$$\frac1{|w|}\Sh(w)=\frac1{|w|}\sum_{a\in A} N_a(w)\Sh(a)\rightarrow \sum_{a\in A} \mu_\varphi(a)\Sh(a)=\Delta_\varphi(\Sh) {\:\rm as}\: |w|\rightarrow\infty.$$

Thus for long words $w$ the images $\Sh(w)$ will be concentrated around the line in the direction of the drift of $\Sh$, and so the complement of $\Sh(\Lan_\varphi)$ will have  infinite cardinality if the drift is not $(0,0)$.
 \hfill$\Box$

\medskip

\qquad Can we say something about the Frobenius problem for homomorphic images of morphic languages of an embedding with drift $(0,0)$?
We shall give an infinite family of morphic languages $\Lan_\theta$ on an alphabet $A=\{a,b,c,d\}$ of four letters where for the homomorphism $\Sh^\oplus$ given by
$$ \Sh^\oplus(a)=(1,0),\; \Sh^\oplus(b)=(0,1),\;\Sh^\oplus(c)=(-1,0),\; \Sh^\oplus(d)=(0,-1)   $$
the homomorphic embedding is the whole $ \Z\times\Z$---and thus the complement is empty.
\qquad We shall make use of the paperfolding morphisms introduced in \cite{Dek-paper}. Let $\sigma$ be the rotation morphism on the alphabet  $\{a,b,c,d\}$ given by
$\sigma(a)=b,\; \sigma(b)=c,\; \sigma(c)=d,\; \sigma(d)=a,$
and let $\tau$ be the anti-morphism given by $\tau(w_1\dots w_n)=w_n\dots w_1$.

A morphism $\theta$ on $\{a,b,c,d\}$ is called a \emph{paperfolding} morphism if

\quad 1) $\sigma \tau \theta=\theta$,

\quad 2) Letters from $\{a,c\}$ alternate\footnote{This corrects an omission in \cite[Definition 1]{Dek-paper}.} with letters from $\{b,d\}$ in $\theta(a)$.

 A paperfolding morphism is called \emph{symmetric} if $\sigma \theta=\theta$.
It is clear that this happens if and only if the word $\theta(a)$ is a palindrome.

\qquad  Let $G$ be a (semi-) group with operation $+$ and unit $e$. In general an infinite word $x=(x_n)$ over an alphabet $A$ and a homomorphism $\Sh:A^*\rightarrow G$ generate a \emph{walk} $Z=(Z_n)_{n\ge 0}$ by (cf.~\cite{Dek-Marches})
$$Z_0=e, \qquad Z_{n+1}=Z_{n}+\Sh(x_n)=\Sh(x_0\dots x_n),\:{\rm for}\: n\ge 0.$$

\qquad A paperfolding morphism $\theta$ with $\theta(a)=a...$ is called \emph{perfect} if the four walks generated by the fixed point $x=\theta^\infty(a)$, and its three rotations over $\pi/2, \pi$ and $3\pi/2$ visit every integer point in the plane exactly twice (except the origin, which is visited 4 times).

In \cite{Dek-paper} it is---not explicitly---proved that for any odd integer $N$  that is the sum of two squares there exists a perfect symmetric paperfolding morphism of length $N$. To make the proof explicit, one uses that according to the paragraph at the end of Section 7 in \cite{Dek-paper} there exists a symmetric planefilling and self-avoiding string for each such $N$, and then
one observes that the construction of such a string in the proof of \cite[Theorem 4]{Dek-paper}  always satisfies
 the perfectness criterion given in \cite[Theorem 5]{Dek-paper}.

\medskip

The smallest length is $N=5$, with morphism $\theta$ given by
$$\theta(a)=abcba,\; \theta(b)=bcdcb, \; \theta(c)=cdadc, \; \theta(d)=dabad.$$

\vspace*{-.4cm}

\begin{figure}[!h]
\centering
\includegraphics[width=6cm]{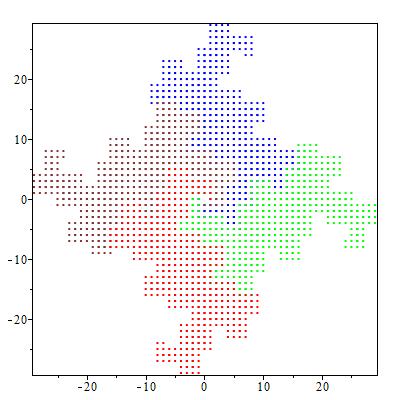}
\caption{\small The four images of the words $\theta^4(a),\dots,\theta^4(d)$ under $\Sh^\oplus$, where $\theta$ is the perfect symmetric 5-folding morphism. The origin is not covered, but it is the image of the word $abcd\in \Lan_\theta$.}\label{fold5}
\end{figure}

\vspace*{.4cm}

 \begin{proposition}\label{prop:drift} Let $\Lan_\theta$ be the language generated by a perfect symmetric paperfolding morphism $\theta$.
 Then  $\Sh^\oplus(\Lan_\theta)= \Z\times \Z$.
\end{proposition}

\emph{Proof:} This follows directly from Theorem 5 in \cite{Dek-paper}, using the observation above.  \hfill$\Box$


%
%
%
%


\end{document}